\documentstyle[11pt,amstex]{amsart}
\newsymbol \varsubsetneq 2320
\newcommand{\forces}
{\mathrel {{\vrule height 6.9pt depth -0.1pt}\! \vdash }}
\newcommand{\rest}
{{\hbox{$\,|\grave{}\,$}}}
\renewcommand{\dot}{\overset{\boldsymbol .}}
\renewcommand{\S}{{\cal S}} 
\renewcommand{\r}{{\Bbb R}}

\newcommand{\comp}{\circ} 
\newcommand{\can}{2^{\textstyle \omega}} 
\newcommand{\baire}{\omega^{\textstyle \omega}} 
 
\newcommand{\fsuo}{[\omega]^{\textstyle <\!\omega}} 
\newcommand{\dom}{
\operatorname {dom}} 
\newcommand{\rng}{ 
\operatorname {rng}}
\newcommand{\nor}{ 
\operatorname {\bold {nor}}\/} 
\newcommand{\A}{{\cal A}}  
\newcommand{\B}{{\cal B}}  
\newcommand{\C}{{\cal C}}
\newcommand{\D}{{\cal D}}

\newcommand{\K}{{\cal K}}
\newcommand{\M}{{\cal M}}  
\newcommand{\N}{{\cal N}}  
\newcommand{\p}{{\Bbb P}}

\newcommand{\q}{{\Bbb Q}}

\newcommand{\Constr}{\noindent{\sc Construction:}\hspace{0.15in}}

\newcommand{\bdom}{\operatorname {\bold {dom}}}
\newcommand{\val}{\operatorname {\bold {val}}}
\newcommand{\sval}{\operatorname {\bold {sval}}}
\newcommand{\dis}{\operatorname {\bold {dis}}}
\newcommand{\pos}{
\operatorname {pos}}
\newcommand{\SCR}{
\operatorname {SCR}}
\newcommand{\bH}{{\bold H}}
\newcommand{\POS}{
\operatorname {POS}}
\newcommand{\qplus}{{\q^+_\infty(K,\Sigma,\Sigma^\bot)}}
\newcommand{\qzero}{{\q_\emptyset(K,\Sigma,\Sigma^\bot)}}
\newcommand{\id}{{\Bbb I}_\infty(K,\Sigma,\Sigma^\bot)}
\newcommand{\trojka}{(K_{\ref{main}},\Sigma_{\ref{main}},\Sigma^\bot_{
\ref{main}})}
\newcommand{\forsing}{\q^+_\infty\trojka}
\newcommand{\idmain}{{\Bbb I}_\infty(K_{\ref{main}},\Sigma_{\ref{main}},
\Sigma_{\ref{main}}^\bot)}
\newcommand{\Leb}{
\operatorname {Leb}}
\newcommand{\qnew}{{\q^{\bH}_{\bar{n}}}}

\newtheorem{theorem}{Theorem}[section] 
\newtheorem{claim}{Claim}[theorem]
\newtheorem{proposition}[theorem]{Proposition} 
\newtheorem{conclusion}[theorem]{Conclusion}

\theoremstyle{definition}
\newtheorem{kprob}[theorem]{Kunen's Problem} 
\newtheorem{definition}[theorem]{Definition}
\newtheorem{notation}[theorem]{Notation}
\newtheorem{example}[theorem]{Example}
\newtheorem{remark}[theorem]{Remark}

\title[\it Norms on possibilities II: more ccc ideals on $\can$] {
\vspace {3.0cm}\uppercase 
{\Large \bf 
Norms on possibilities II: more ccc ideals on $\can$}} 
\author[A. Ros{\l}anowski]{\uppercase {\bf A. Ros{\l}anowski}}
\author[S. Shelah]{\uppercase {\bf S. Shelah}}

\thanks{The first author thanks the Hebrew University and the Lady Davis
Foundation for the Golda Meir Postdoctoral Fellowship. The research of the
second author was partially supported by ``Basic Research Foundation'' of the
Israel Academy of Sciences and Humanities. Publication 628 of the second
author.}
\subjclass{\ } 

\begin{document}
\ 

\ 

\maketitle 

\bigskip
\bigskip
\bigskip
{\it Abstract.}
We use the method of {\em norms on possibilities} to answer a question of
Kunen and construct a ccc $\sigma$--ideal on $\can$ with various closure
properties and distinct from the ideal of null sets, the ideal of meager sets
and their intersection.

\bigskip
\bigskip
\bigskip
\stepcounter{section}
\subsection*{\quad 0. Introduction}
In the present paper we use the method of {\em norms on possibilities} to
answer a question of Kunen (see \ref{pwa} below) and construct a ccc
$\sigma$--ideal on $\can$ with various closure properties and distinct from
the ideal of null sets, the ideal of meager sets and their intersection. The
method we use is, in a sense, a generalization of the one studied
systematically in \cite{RoSh:470} (the case of creating pairs). However, as
the main desired property of the forcing notion we construct is satisfying the
ccc, we do not use the technology of that paper (where the forcing notions
were naturally proper not--ccc) and our presentation does not require
familiarity with the previous part.  

The following problem was posed by Kunen and for some time stayed open.

\medskip
\begin{kprob}
[\mbox {see \cite[Question 1.2]{Ku84} or \cite[Problem 5.5]{Mi91}}]
\label{pwa}
Do\-es \linebreak 
there exist a ccc $\sigma$--ideal ${\cal I}$ of subsets of $\can$ 
which:
\begin{quotation}
\noindent has a Borel basis, is index--invariant, translation--invariant and
is neither the meager ideal, nor the null ideal nor their intersection?
\end{quotation}
\end{kprob}

\medskip
We give a positive answer to this problem and, in fact, we construct a large
family of ideals with the desired properties. Naturally, these ideals are
obtained from ccc forcing notions. However, the construction presented in
Section 4 is purely combinatorial and no forcing techniques are needed for it.

Let us recall the definitions of the properties required from ${\cal I}$ in
\ref{pwa}. 

\medskip
\begin{definition}
\label{idedef}
\ 

\begin{enumerate}
\item If $({\cal X},+)$ is a commutative group then the product space ${\cal X
}^{\textstyle \omega}$ is equipped with the product group structure. The
Cantor space $\can$ carries the structure of Polish group generated by the
addition modulo 2 on $2=\{0,1\}$. We equip the Baire space $\baire$ with a
group operation interpreting $\omega$ as the group ${\Bbb Z}$ of integers. 
\item An ideal ${\cal I}$ of subsets of a group $({\cal X},+)$ is {\em
translation--invariant} if 
\[(\forall x\in{\cal X})(\forall A\in {\cal I})(A+x\stackrel {\rm def}{=}
\{a+x:a\in A\}\in {\cal I}).\]
\item An ideal ${\cal I}$ of subsets of the product space ${\cal
X}^{\textstyle \omega}$ is {\em index--invariant} if for every embedding $\pi:
\omega\stackrel{1-1}{\longrightarrow}\omega$ and every set $A\in {\cal I}$ we
have 
\[\pi_*(A)\stackrel {\rm def}{=}\{x\in {\cal X}^\omega: x\comp\pi\in A\}\in
{\cal I}.\]
\item The ideal ${\cal I}$ as above is {\em permutation--invariant} if it
satisfies the demand above when we restrict ourselves to permutations $\pi$ of
$\omega$ only.
\item An ideal ${\cal I}$ on a Polish space ${\cal X}$ has a Borel basis if
every set $A\in {\cal I}$ is contained in a Borel set $B\in {\cal I}$.

\noindent [In this situation we may say that ${\cal I}$ is a Borel ideal.]
\item A Borel $\sigma$-ideal ${\cal I}$ on a Polish space ${\cal X}$ is ccc if
the quotient Boolean algebra $
\text {BOREL}({\cal X})/{\cal I}$ of Borel subsets
of ${\cal X}$ modulo $\cal I$ satisfies the ccc. (Equivalently, there is no
uncountable family of disjoint Borel subsets of ${\cal X}$ which are not in
$\cal I$.)  
\end{enumerate}
\end{definition}

\medskip

There has been some partial answers to \ref{pwa} already. Kechris and Solecki 
\cite{KeSo95} showed that ccc $\sigma$--ideals generated by closed sets
(i.e.~with $\Sigma^0_2$--basis) are essentially like the meager ideal. 
It was shown in \cite{JuRo92} how Souslin ccc forcing notions may produce nice
$\sigma$--ideals on the Baire space $\baire$ (those notes were presented in
\cite[pp 193--203]{BaJu95}). That method provided an answer to \ref{pwa} if
one replaced the demand that the ideal in question is on $\can$ by allowing it
to be on $\baire$. It seems that the approach presented there is not
applicable if we want to stay in the Cantor space. Our method here, though
similar to the one there, is more direct. 

\medskip
\noindent{\bf Notation:} Our notation is rather standard and compatible with
that of classical textbooks on Set Theory (like Jech \cite{J} or
Bartoszy\'nski Judah \cite{BaJu95}). However in forcing we keep the convention
that {\em a stronger condition is the larger one}. 

\medskip
\begin{notation}
\label{notacja}
\

\begin{enumerate}
\item $\r^{{\geq}0}$ stands for the set of non-negative reals. The integer
part of a real $r\in\r^{{\geq}0}$ is denoted by $\lfloor r\rfloor$.
\item For a set $X$,\ \ \ $[X]^{\textstyle{\leq}\omega}$,
$[X]^{\textstyle{<}\omega}$ and ${\cal P}(X)$ will stand for families of
countable, finite and all, respectively, subsets of the set $X$. The family of
$k$-element subsets of $X$ will be denoted by $[X]^{\textstyle k}$. The set of
all finite sequences with values in $X$ is called $X^{\textstyle {<}\omega}$
(so domains of elements of $X^{\textstyle {<}\omega}$ are integers). 
\item The Cantor space $\can$ and the Baire space $\baire$ are the spaces of
all functions from $\omega$ to $2$, $\omega$, respectively, equipped with
natural (Polish) topology. 
\item For a forcing notion $\p$, $\Gamma_\p$ stands for the canonical
$\p$--name for the generic filter in $\p$. With this one exception, all
$\p$--names for objects in the extension via $\p$ will be denoted with a dot
above (e.g. $\dot{\tau}$, $\dot{X}$).
\end{enumerate}
\end{notation}

\medskip
\noindent {\bf Basic Notation:} In this paper $\bH$ will stand for a function
with domain $\omega$ such that $(\forall m\in\omega)(|\bH(m)|\geq 2)$. We
usually assume that $0\in \bH(m)$ (for all $m\in\omega$). Moreover, we assume
that at least $\bH\in {\cal H}(\aleph_1)$ (the family of hereditarily
countable sets) or, what is more natural, even $\bH(i)\in\omega\cup\{\omega
\}$ (for $i\in\omega$).

\bigskip
\stepcounter{section}
\subsection*{\quad 1. Semi--creating triples}
In \cite{RoSh:470} we explored a general method of building forcing notions
using {\em norms on possibilities}. We studied weak creating pairs and their
two specific cases: creating pairs and tree--creating pairs. For our
applications here, we have to modify the general schema introducing an
additional operation. In our presentation we will not refer the reader to
\cite{RoSh:470}, but familiarity with that paper may be of some help in
getting a better picture of the method.

The definition \ref{semi} of semi--creatures and semi--creating triples is not
covered by the general case as presented in \cite[1.1]{RoSh:470}. But one
should notice a close relation of it to the case discussed in
\cite[1.2]{RoSh:470} under an additional demand that the considered creating
pair is forgetful (see \cite[1.2.5]{RoSh:470}). 

\medskip
\begin{definition}
\label{semi}
Let $\bH:\omega\longrightarrow {\cal H}(\aleph_1)$.
\begin{enumerate}
\item A quadruple $t=(\bdom,\sval,\nor,\dis)$ is {\em a semi--creature for
$\bH$} if
\begin{itemize}
\item $\bdom\in \fsuo\setminus\{\emptyset\}$,
\item $\sval\subseteq\prod\limits_{i\in\bdom}\bH(i)$ is non-empty,
\item $\nor\in \r^{\geq 0}\cup\{\infty\}$,
\item $\dis\in {\cal H}(\aleph_1)$, and
\item $\sval=\prod\limits_{i\in\bdom}\bH(i)$\quad if and only if\quad
$\nor=\infty$.
\end{itemize}
The family of all semi--creatures for $\bH$ will be denoted by $\SCR[\bH]$.
\item In the above definition we write $\bdom=\bdom[t]$, $\sval=\sval[t]$,
$\nor=\nor[t]$ and $\dis=\dis[t]$.
\item Suppose that $K\subseteq \SCR[\bH]$. A function $\Sigma:[K]^{\textstyle
<\omega}\longrightarrow {\cal P}(K)$ is {\em a semi--composition operation on
$K$} if for each $\S\in [K]^{\textstyle<\omega}$ and $t\in K$:
\begin{enumerate}
\item[(a)] if  for each $s\in\S$ we have $s\in\Sigma(\S_s)$, $S_s\in [K]^{
\textstyle<\omega}$, then $\Sigma(\S)\subseteq \Sigma(\bigcup\limits_{s\in\S}
\S_s)$,
\item[(b)] $t\in\Sigma(t)$, $\Sigma(\emptyset)=\emptyset$, and
\item[(c)] if $t\in\Sigma(\S)$ (so $\Sigma(\S)\neq\emptyset$) then
\begin{enumerate}
\item[($\alpha$)] $\bdom[t]=\bigcup\limits_{s\in\S}\bdom[s]$,
\item[($\beta$)]  $v\in\sval[t]\ 
\&\ 
s\in\S\quad\Rightarrow\quad v \rest
\bdom[s]\in\sval[s]$, and 
\item[($\gamma$)] $s_1,s_2\in\S\ 
\&\ 
s_1\neq s_2\quad\Rightarrow\quad \bdom[
s_1]\cap\bdom[s_2]=\emptyset$.
\end{enumerate}
\end{enumerate}
\item A mapping $\Sigma^\bot:K\longrightarrow [K]^{\textstyle{<}\omega}
\setminus\{\emptyset\}$ is called {\em a semi--decomposition operation on $K$} 
if for each $t\in K$: 
\begin{enumerate}
\item[(a)$^\bot$]  if $\S=\{s_0,\ldots,s_k\}\in\Sigma^\bot(t)$ and
$\S_i\in\Sigma^\bot(s_i)$ (for $i\leq k$) then $\S_0\cup\ldots\cup\S_k\in
\Sigma^\bot(t)$, 
\item[(b)$^\bot$] $\{t\}\in\Sigma^\bot(t)$,
\item[(c)$^\bot$]  if $\S\in\Sigma^\bot(t)$ then 
\begin{enumerate}
\item[$(\alpha)^\bot$] $\bdom[t]=\bigcup\limits_{s\in\S}\bdom[s]$,
\item[$(\beta)^\bot$]  $\{v\in\prod\limits_{i\in\bdom[t]}\bH(i):(\forall
s\in\S)(v \rest \bdom[s]\in\sval[s])\}\subseteq \sval[t]$. 
\end{enumerate}
\end{enumerate}
\item If $\Sigma$ is a semi--composition operation on $K\subseteq\SCR[\bH]$
and $\Sigma^\bot$ is a semi--decomposition on $K$ then $(K,\Sigma,
\Sigma^\bot)$ is called {\em a semi--creating triple for $\bH$}.  
\item If we omit $\bH$ this means that either $\bH$ should be clear from the
context or we mean {\em for some} $\bH$. 
\end{enumerate}
\end{definition}

\medskip
\begin{remark}
\ 

\begin{enumerate}
\item In the definition of semi--creatures above, $\bdom$ stands for domain,
$\sval$ for semi--values, $\nor$ is for norm and $\dis$ for distinguish. The
last plays a role of an additional parameter and it may be forgotten sometimes
(compare \cite[1.1.2]{RoSh:470}). One should notice that the difference with
\cite[1.1.1, 1.2.1]{RoSh:470} is in $\bdom$ (which in case of creatures of
\cite[1.2]{RoSh:470} is always an interval). This additional freedom has some 
price: we put a slightly more restrictive demands on $\val$, so we have
$\sval$ here. We could have more direct correspondence between $\val$ of 
\cite[1.1]{RoSh:470} and $\sval$ here, but that would complicate notation only.
\item Note that in \ref{semi}(3) we allow $\Sigma(\S)=\emptyset$ even if
$\S\neq\emptyset$. In applications we will say what are the values of
$\Sigma(\S)$ only if it is a non-empty set; so not defining $\Sigma(\S)$ means
that the value is $\emptyset$. The demand \ref{semi}(3c) appears to simplify
\ref{forcing} below only; we could have used there $\pos$, defined like in
\cite[1.1.6]{RoSh:470}. 
\item The main innovation here is the additional operation $\Sigma^\bot$,
which will play a crucial role in getting the ccc. If it is trivial
(i.e.~$\Sigma^\bot(t)=\{\{t\}\}$) then we are almost in the case of creating
pairs of \cite[1.2]{RoSh:470}.
\end{enumerate}
\end{remark}

\medskip
\begin{definition}
\label{forcing}
Assume $(K,\Sigma,\Sigma^\bot)$ is a semi--creating triple for $\bH$. We
define forcing notions $\qzero$ and $\qplus$ as follows.

\medskip
\noindent{\bf 1)\quad Conditions} in $\qzero$ are sequences
$(w,t_0,t_1,t_2,\ldots)$ such that 
\begin{enumerate}
\item[(a)] $w\in\prod\limits_{i\in\dom(w)}\bH(i)$ is a finite function,
\item[(b)] each $t_i$ belongs to $K$  and $\omega=\dom(w)\cup
\bigcup\limits_{i\in\omega}\bdom[t_i]$ is a partition of $\omega$ (so
$i<j<\omega$ implies $\dom(w)\cap\bdom[t_i]=\bdom[t_i]\cap\bdom[t_j]=
\emptyset$).
\end{enumerate}

\noindent{\bf The relation $\leq_\emptyset$} on $\qzero$ is given by:
\quad $(w_1,t^1_0,t^1_1,t^1_2,\ldots)\leq_\emptyset
(w_2,t^2_0,t^2_1,t^2_2,\ldots)$
\quad if and only if\quad $(w_2,t^2_0,t^2_1,t^2_2,\ldots)$ can be obtained
from $(w_1,t^1_0,t^1_1,t^1_2,\ldots)$ by applying finitely many times the
following operations (in the description of the operations we say what are
their legal results for a condition $(w,t_0,t_1,t_2,\ldots)\in\qzero$).  

{\em Deciding the value} for $(w,t_0,t_1,t_2,\ldots)$:

\noindent a legal result is a condition $(w^*,t^*_0,t^*_1,t^*_2,\ldots)\in
\qzero$ such that for some finite $A\subseteq\omega$ (possibly empty) we have\\
$w\subseteq w^*$,\quad $\dom(w^*)=\dom(w)\cup\bigcup\limits_{i\in A}
\bdom[t_i]$,\quad $w^* \rest \bdom[t_i]\in\sval[t_i]$ (for $i\in A$)\quad 
and $\{t^*_0,t^*_1,\ldots\}=\{t_i: i\in\omega\setminus A\}$. 

{\em Applying $\Sigma$} to $(w,t_0,t_1,t_2,\ldots)$:

\noindent a legal result is a condition $(w,t^*_0,t^*_1,t^*_2,\ldots)\in
\qzero$ such that for some disjoint sets $A_0,A_1,A_2,\ldots\in\fsuo$ we have
\[t^*_i\in\Sigma(t_j:j\in A_i)\quad\mbox{ for each }i\in\omega.\]

{\em Applying $\Sigma^\bot$} to $(w,t_0,t_1,t_2,\ldots)$:

\noindent a legal result is a condition $(w,t^*_0,t^*_1,t^*_2,\ldots)\in
\qzero$ such that for some disjoint non-empty finite sets $A_0,A_1,A_2,\ldots
\subseteq\omega$ we have $\{t^*_j: j\in A_i\}\in\Sigma^\bot(t_i)$.
\medskip

\noindent{\bf 2)}\quad If $p=(w,t_0,t_1,t_2,\ldots)\in\qzero$ then we let
$w^p=w$, $t^p_i=t_i$. 
\medskip

\noindent{\bf 3)}\quad Finally we define $\qplus$ as the collection of all
$p\in\qzero$ such that  
\begin{enumerate}
\item[(c)] $(\forall i\in\omega)(\nor[t^p_i]\neq\infty)$\quad and\quad 
$\lim\limits_{i\to\infty}\nor[t^p_i]=\infty$,
\end{enumerate}
and the relation $\leq$ on $\qplus$ is the restriction of $\leq_\emptyset$.
\end{definition}

\medskip
\begin{proposition}
If $(K,\Sigma)$ is a semi--creating pair then $\qzero,\qplus$ are forcing
notions (i.e.~the relation $\leq_\emptyset$ is transitive). 
\end{proposition}

\medskip
\begin{remark}
\ 

\begin{enumerate}
\item Like in \cite{RoSh:470}, one may consider various variants of the demand
\ref{forcing}(3)(c).
\item Note that in a condition in $\qplus$ we allow only those $t_i\in K$ for
which $\nor[t_i]\neq\infty$. We could restrict $K$ to semi--creatures with
finite norm, but presence of $t$ which give no restrictions will make some
definitions simpler. 
\item One should notice a close relation of $\qplus$ here to forcing notions
of type $\q^*_\infty(K,\Sigma)$ in \cite[1.2]{RoSh:470}. What occurs in
applications here, is that (in interesting cases) the forcing notion $\qplus$
is equivalent to some $\q^*_\infty(K',\Sigma')$.
\end{enumerate}
\end{remark}

\medskip
As in the present paper we are interested mainly in $\sigma$--ideals on the
real line, let us show how our forcing notions introduce ideals on
$\prod\limits_{i\in\omega}\bH(i)$. Later we will look more closely at $\qplus$
as a forcing notion. 

\medskip
\begin{definition}
\label{ideals}
Assume $(K,\Sigma,\Sigma^\bot)$ is a semi--creating triple for $\bH$.
\begin{enumerate}
\item For a condition $p\in\qplus$ we define its {\em total possibilities} as 
\[\POS(p)=\{x\in\prod_{i\in\omega}\bH(i): w^p\subseteq x\ \&\ (\forall j\in
\omega)(x \rest
\bdom[t^p_j]\in\sval[t^p_j])\}.\]
\item Let $\dot{W}$ be a $\qplus$--name such that
\[ \forces
_{\qplus}\dot{W}=\bigcup\{w^p: p\in\Gamma_{\qplus}\}\]
(compare to \cite[1.1.13]{RoSh:470}).
\item Let $\id$ be the ideal of subsets of $\prod\limits_{i\in\omega}\bH(i)$
generated by those Borel sets $B\subseteq\prod\limits_{i\in\omega}\bH(i)$ for
which $ \forces
_{\qplus}$``$\dot{W}\notin B$''. 

\noindent [Remember that $\bH(i)$ is countable (for each $i\in\omega$); so
considering the product topology on $\prod\limits_{i\in\omega}\bH(i)$ (with each
$\bH(i)$ discrete) we get a Polish space.] 
\end{enumerate}
\end{definition}

\medskip
\begin{proposition}
\label{idccc}
Assume $(K,\Sigma,\Sigma^+)$ is a semi--creating triple for $\bH$.
\begin{enumerate}
\item If $p,q\in\qplus$, $p\leq q$ then $\POS(q)\subseteq \POS(p)$, $\POS(p)$
is a non-empty closed subset of $\prod\limits_{i\in\omega}\bH(i)$ and
$p \forces $``$\dot{W}\in\POS(p)$''.
\item $\id$ is a $\sigma$--ideal, $\prod\limits_{i\in\omega}\bH(i)\notin\id$
(in fact, \linebreak $\POS(p)\notin\id$ for $p\in\qplus$).
\item If the forcing notion $\qplus$ satisfies the ccc then the ideal $\id$
satisfies the ccc. 
\end{enumerate}
\end{proposition}

\medskip
To make sure that the ideal $\id$ is invariant we have to assume natural
invariance properties of the semi--creating triple.

\medskip
\begin{definition}
\label{invdef}
Assume that $\bH$ is a constant function, say $\bH(i)={\cal X}$, and we have a
commutative group operation $+$ on ${\cal X}$. Let $(K,\Sigma,\Sigma^\bot)$ be
a semi--creating triple for $\bH$. We say that 
\begin{enumerate}
\item $(K,\Sigma,\Sigma^\bot)$ is {\em directly $+$--invariant} if for each
$t\in K$ and $v\in {\cal X}^{\textstyle\bdom[t]}$ there is a unique
semi--creature $s\in K$ (called $t+v$) such that
\[\bdom[s]=\bdom[t],\quad \nor[s]=\nor[t],\]
\[\sval[s]=\sval[t]+v=\{w+v: w\in\sval[t]\},\]
and $\Sigma^\bot(t+v)=\{\{s+(v \rest 
\bdom[s]): s\in\S\}:\S\in\Sigma^\bot(t)\}$
(for $t\in K$, $v\in {\cal X}^{\textstyle \bdom[t]}$),\\
and $\Sigma(t+(v \rest
\bdom[t]):t\in\S)=\{s+v: s\in\Sigma(\S)\}$ (for $\S\in
[K]^{\textstyle{<}\omega}\setminus\{\emptyset\}$ and $v:\bigcup\limits_{s\in
\S} \bdom[s]\longrightarrow{\cal X}$),
\item $(K,\Sigma)$ is {\em directly permutation--invariant} if for each $t\in
K$ and an embedding $\pi:X\stackrel{1-1}{\longrightarrow}\omega$ such that
$\bdom[t]\subseteq X\subseteq\omega$ there is a unique semi--creature $s\in K$
(called $\pi(t)$) such that 
\[\bdom[s]=\pi[\bdom[t]],\quad \nor[s]=\nor[t],\]  
\[\sval[s]=\{w\comp\pi^{-1}: w\in \sval[t]\},\]
and $\Sigma^\bot(\pi(t))=\{\{\pi(s): s\in\S\}:\S\in\Sigma^\bot(t)\}$ (for $t
\in K$, $\pi:\omega\stackrel{1-1}{\longrightarrow}\omega$), \\
and $\Sigma(\pi(t):t\in\S)=\{\pi(s):s\in\Sigma(\S)\}$ (for $\S\in [K]^{
\textstyle {<}\omega}$, $\pi:\omega\stackrel{1-1}{\longrightarrow}\omega$).
\end{enumerate}
\end{definition}

\medskip
\begin{proposition}
\label{invthm}
Assume that $\bH$ is a constant function (and $\bH(i)={\cal X}$ for
$i\in\omega$), and we have a commutative group operation $+$ on ${\cal X}$ (so
then $\prod\limits_{i\in\omega}\bH(i)={\cal X}^{\textstyle \omega}$ becomes a
commutative group too). Let $(K,\Sigma,\Sigma^\bot)$ be a semi--creating
triple for $\bH$. 
\begin{enumerate}
\item If $(K,\Sigma,\Sigma^\bot)$ is directly $+$--invariant then the ideal
$\id$ is translation invariant (in the product group ${\cal X}^{\textstyle
\omega}$).
\item If $(K,\Sigma,\Sigma^\bot)$ is directly permutation--invariant then
$\id$ is permutation invariant (see \ref{idedef}(4)). 
\end{enumerate}
\end{proposition}

\begin{pf} 
\quad 1)\ \ \ Note that if $v\in {\cal X}^{\textstyle\omega}$ then the
mapping 
\[(w,t_0,t_1,\ldots)\mapsto(w+(v \rest
\dom(w)), t_0+(v \rest
\bdom[t_0]),t_1
\rest
(\bdom[t_1]),\ldots)\]
is an automorphism of the forcing notion $\qplus$.

\noindent 2) Similarly. 
\end{pf}

\medskip
In general it is not clear if the ideal $\id$ contains any non-empty set. One
can easily formulate a demand ensuring this (like \cite[3.2.7]{RoSh:470}), but
there is no need for us to deal with it, as anyway we want to finish with a
ccc ideal. 

\bigskip
\stepcounter{section}
\subsection*{\quad 2. Getting ccc}
In this section we show how we can make sure that the forcing notion $\qplus$
satisfies the ccc. The method for this (and the required properties) are quite
simple and they will allow us to conclude more properties of the ideal
$\id$. The main difficulty will be to construct an example meeting all the
appropriate demands (and this will be done in the next section). 

\medskip
\begin{definition}
\label{cccdef}
Assume that $(K,\Sigma,\Sigma^\bot)$ is a semi--creating triple for $\bH$.
\begin{enumerate}
\item We say that {\em $(K,\Sigma,\Sigma^\bot)$ is linked} if for each
$t_0,t_1\in K$ such that  
\[\nor[t_0],\nor[t_1]>1\quad\mbox{ and }\quad\bdom[t_0]=\bdom[t_1]\]
there is $s\in \Sigma(t_0)\cap\Sigma(t_1)$ with 
\[\nor[s]\geq\min\{\nor[t_0],\nor[t_1]\}-1.\]
\item The triple $(K,\Sigma,\Sigma^\bot)$ is called {\em semi--gluing} if for
each $t_0,\ldots,t_n\in K$ such that $k<\ell\leq n\quad\Rightarrow\quad
\bdom[t_k]\cap\bdom[t_\ell]=\emptyset$ there is $s\in \Sigma(t_0,\ldots,t_n)$
with 
\[\nor[s]\geq\min\{\nor[t_k]: k\leq n\}-1.\]
\item We say that {\em $(K,\Sigma,\Sigma^\bot)$ has the cutting property} if
for every $t\in K$ with $\nor[t]>1$ and for each non-empty $z
\varsubsetneq
\bdom[
t]$ there are $s_0,s_1\in K$ such that
\begin{enumerate}
\item[$(\alpha)$] $\bdom[s_0]=z$, $\bdom[s_1]=\bdom[t]\setminus z$,
\item[$(\beta)$]  $\nor[s_\ell]\geq \nor[t]-1$ (for $\ell=0,1$) and
\item[$(\gamma)$] $\{s_0,s_1\}\in\Sigma^\bot(t)$.
\end{enumerate}
\end{enumerate}
\end{definition}

\medskip
\begin{remark}
\ 

\begin{enumerate}
\item Semi--gluing triples $(K,\Sigma,\Sigma^\bot)$ correspond to gluing
creating pairs (as defined in \cite[2.1.7(2)]{RoSh:470}).
\item Note that in \ref{cccdef}(3) we do not require that $\nor[s_\ell]\neq
\infty$. However, if $t\in K$ satisfies $\nor[t]\neq\infty$ and $s_0,s_1$ are
as in \ref{cccdef}(3) then necessarily at least one of the $s_\ell$'s has to
have the same property. 
\end{enumerate}
\end{remark}

\medskip
\begin{definition}
\label{linked}
A forcing notion $\q$ is {\em $\sigma$-$*$--linked} if for every $n\in\omega$
there is a partition $\langle A_i: i\in\omega\rangle$ of $\q$ such that
\[\mbox{if }\quad q_0,\ldots,q_n\in A_i\quad\mbox{ then }\quad (\exists q\in
\q)(q_0\leq q\ \&\ \ldots\ \&\ q_n\leq q).\]
\end{definition}

\medskip
\begin{theorem}
\label{cccthm}
Assume that $(K,\Sigma,\Sigma^\bot)$ is a linked semi--gluing semi--creating
triple with the cutting property. Then the forcing notion $\qplus$ is
$\sigma$-$*$--linked.
\end{theorem}

\begin{pf}
\quad  Fix $n\in\omega$. For a finite function $w\in\prod\limits_{i\in
\dom(w)}\bH(i)$ let 
\[A_w\stackrel {\rm 
def}{=}\{p\in\qplus: w^p=w\ \&\ (\forall i\in\omega)(\nor[
t^p_i]>n+5)\}.\]
Clearly $\bigcup\{A_w: w\in\prod\limits_{i\in\dom(w)}\bH(i),
\dom(w)\in\fsuo\}$ is a dense subset of $\qplus$, so it is enough to show the
following claim. 

\begin{claim}
\label{cl1}
Let $w\in\prod\limits_{i\in\dom(w)}\bH(i)$, $\dom(w)\in\fsuo$ and
$p_0,\ldots,p_n\in A_w$. Then the conditions $p_0,\ldots,p_n$ have a common
upper bound (in the forcing notion $\qplus$).
\end{claim}

\noindent{\em Proof of the claim:}\quad Suppose $p_0,\ldots,p_n\in A_w$ (so in
particular $w=w^{p_0}=\ldots=w^{p_n}$). Choose an increasing sequence $\langle
m_i: i<\omega\rangle\subseteq\omega$ such that $\dom(w)\subseteq m_0$ and for
every $\ell\leq n$ and each $i\in\omega$:
\[\begin{array}{l}
(\exists j\in\omega)(\bdom[t^{p_\ell}_j]\subseteq m_0),\quad
(\exists j\in\omega)(\bdom[t^{p_\ell}_j]\subseteq [m_i,m_{i+1}))\quad\mbox{
and}\\
(\forall j\in\omega)(\bdom[t^{p_\ell}_j]\cap m_i\neq\emptyset\ \ \Rightarrow\
\ \bdom[t^{p_\ell}_j]\subseteq m_{i+1})\quad\mbox{ and}\\
(\forall j\in\omega)(\bdom[t^{p_\ell}_j]\setminus m_i\neq\emptyset\ \
\Rightarrow\ \ \nor[t^{p_\ell}_j]>n+5+i).
  \end{array}\] 
Fix $\ell\leq n$ for a moment.\\
For each $j\in\omega$ and $i\in\omega$ such that $\bdom[t^{p_\ell}_j]\cap m_i
\neq \emptyset$ and $\bdom[t^{p_\ell}_j]\cap [m_i,m_{i+1})\neq\emptyset$ use
\ref{cccdef}(3) to choose $s^{\ell,0}_j,s^{\ell,1}_j\in K$ such that
\begin{itemize}
\item $\bdom[s^{\ell,0}_j]=\bdom[t^{p_\ell}_j]\cap m_i$, $\bdom[s^{\ell,1}_j]=
\bdom[t^{p_\ell}_j]\cap [m_i,m_{i+1})$, 
\item $\nor[s^{\ell,0}_j]\geq\nor[t^{p_\ell}_j]-1$, $\nor[s^{\ell,1}_j]\geq
\nor[t^{p_\ell}_j]-1$,
\item $\{s^{\ell,0}_j,s^{\ell,1}_j\}\in\Sigma^\bot(t^{p_\ell}_j)$.
\end{itemize}
Let $\S^\ell_0$ consist of all $t^{p_\ell}_j$ (for $j<\omega$) such that
$\bdom[t^{p_\ell}_j]\subseteq m_0$ and all $s^{\ell,0}_j$ such that
\[\bdom[t^{p_\ell}_j]\cap m_0\neq\emptyset\neq\bdom[t^{p_\ell}_j]\cap [m_0,
m_1).\] 
It should be clear that elements of $\S^\ell_0$ have disjoint domains and
$\bigcup\limits_{s\in\S^\ell_0}\bdom[s]=m_0\setminus\dom(w)$. Use
\ref{cccdef}(2) to find $r^\ell_0\in\Sigma(\S^\ell_0)$ such that
$\nor[r^\ell_0]>n+3$ (re\-mem\-ber the definition of $A_w$). Note that
necessarily $\nor[r^\ell_0]\neq\infty$ as there is $j$ such that
$\bdom[t^{p_\ell}_j]\subseteq m_0$ (remember \ref{semi}(3) and
\ref{forcing}(3)). Similarly, for each $i>0$ we take $\S^\ell_i$ to be the
collection of all $t^{p_\ell}_j$ such that $\bdom[t^{p_\ell}_j]\subseteq
[m_{i-1},m_i)$ and all $s^{\ell,1}_j$ such that
\[\bdom[t^{p_\ell}_j]\cap m_{i-1}\neq\emptyset\neq\bdom[t^{p_\ell}_j]\cap
[m_{i-1},m_i),\] 
and all $s^{\ell,0}_j$ such that 
\[\bdom[t^{p_\ell}_j]\cap m_i\neq\emptyset\neq\bdom[t^{p_\ell}_j]\cap [m_i,
m_{i+1}).\]
Now apply \ref{cccdef}(2) to get $r^\ell_i\in\Sigma(\S^\ell_i)$ such that
$\nor[r^\ell_i]>n+2+i$ (remember the choice of the sequence $\langle m_i:
i<\omega\rangle$; note that, like before, $\nor[r^\ell_i]\neq\infty$). 

It should be clear that $(w,r^\ell_0,r^\ell_1,r^\ell_2,\ldots)$ is a condition
in $\qplus$ stronger than $p_\ell$. Moreover, for each $\ell\leq n$, 
\[\bdom[r^\ell_0]=m_0\setminus\dom(w)\quad\mbox{ and }\quad\bdom[r^\ell_{i+1}]
= [m_i,m_{i+1}).\]
By \ref{cccdef}, for each $i\in\omega$ we find $s^*_i\in K$ such that
\[s^*_i\in\Sigma(r^0_i)\cap\ldots\cap\Sigma(r^n_i)\quad\mbox{ and }\quad
\nor[s^*_i]>i+2.\]
Look at $(w,s^*_0,s^*_1,s^*_2,\ldots)$ --- it is a condition in $\qplus$
stronger than all $p_0,\ldots,p_n$. The claim and the theorem are proved. 
\end{pf}

\bigskip
\stepcounter{section}
\subsection*{\quad 3. The example}
In this section we construct a semi--creating triple with all nice properties
defined and used in the previous section.
 
\medskip
\begin{example}
\label{main}
Let $({\cal X},+)$ be a commutative group, $2\leq |{\cal X}|\leq\omega_0$,
and let $\bH(i)={\cal X}$. Then there is a semi--creating triple $\trojka$
for $\bH$ such that 
\begin{enumerate}
\item the forcing notion $\forsing$ is not trivial,
\item $\trojka$ is directly +--invariant, directly permutation--in\-va\-riant,
\item it is linked and semi--gluing,
\item it has the cutting property.
\end{enumerate}
\end{example}

\Constr Let $\K$ consist of all pairs $(z,\Delta)$ such that
\begin{enumerate}
\item[(a)] $z$ is a non-empty finite subset of $\omega$, and
\item[(b)] $\Delta$ is a non-empty set of non-empty partial functions such
that \linebreak
$\dom(\eta)\subseteq z$ and $\rng(\eta)\subseteq {\cal X}$ for $\eta\in
\Delta$. 
\end{enumerate}
For $(z,\Delta)\in\K$ we define
\[\begin{array}{ll}
v(z,\Delta)\stackrel {\rm 
def}{=}&\{x\in {\cal X}^{\textstyle z}: \neg(\exists
\eta\in\Delta)(\eta\subseteq x)\}\\
n(z,\Delta)\stackrel {\rm 
def}{=}&\max\{k\in\omega:\mbox{for every }\Delta'
\subseteq\Delta\mbox{ there is }\Delta''\subseteq\Delta'\mbox{ such that}\\
\ &\qquad\qquad\mbox{elements of $\Delta''$ have pairwise disjoint domains}\\
\ &\qquad\qquad\mbox{and }\ |\bigcup\limits_{\eta\in\Delta''}\dom(\eta)|\geq
k\cdot |\Delta'|\}. 
  \end{array}\]
Note that the set in the definition of $n(z,\Delta)$ contains 0. Elements of
$\K$ and the two functions $n,v$ are the main ingredients of our
construction. Before we define $\trojka$ let us show some properties of
$(\K,n,v)$. 

\begin{claim}
\label{cl2}
\ 

\begin{enumerate}
\item If $(z,\Delta)\in\K$, $n(z,\Delta)>0$ then $v(z,\Delta)\neq\emptyset$.
\item For each non-empty $z\in\fsuo$ there is $\Delta$ such that
$(z,\Delta)\in\K$ and $n(z,\Delta)=|z|$.
\end{enumerate}
\end{claim}

\noindent{\em Proof of the claim:}\quad 1)\ \ \  Suppose $n(z,\Delta)>0$. 
Then we know
that for each $\Delta'\subseteq\Delta$ there is $\Delta''\subseteq\Delta'$
such that
\[\eta_0,\eta_1\in\Delta''\ \ \&\ \ \eta_0\neq\eta_1\quad\Rightarrow\quad
\dom(\eta_0)\cap\dom(\eta_1)=\emptyset\] 
and $|\bigcup\limits_{\eta\in\Delta''}\dom(\eta)|\geq|\Delta'|$. Thus we may
use the marriage theorem of Hall (see \cite{Ha35}) and choose a system of
distinct representatives of $\{\dom(\eta):\eta\in\Delta\}$. So for
$\eta\in\Delta$ we have $x_\eta\in\dom(\eta)$ such that
\[\eta_0,\eta_1\in\Delta\ \ \&\ \ \eta_0\ne\eta_1\quad\Rightarrow\quad
x_{\eta_0}\neq x_{\eta_1}.\]
Let $w\in {\cal X}^{\textstyle z}$ be such that $w(x_\eta)\in {\cal X}
\setminus\{\eta(x_\eta)\}$ for $\eta\in\Delta$. Clearly $w\in v(z,\Delta)$. 

\medskip
\noindent 2)\ \ \ Take $\Delta=\{{\bf 0}_z\}$, where ${\bf 0}_z$ is a function
on $z$ with constant value $0$.

\begin{claim}
\label{cl3}
\ 

\begin{enumerate}
\item Suppose that $(z,\Delta)\in\K$, $\emptyset\neq z^*\subseteq z$ and 
\[\Delta^*\stackrel {\rm 
def}{=}\{\eta \rest 
z^*: \eta\in\Delta\ \&\ |\dom(\eta)
\cap z^*|\geq\frac{1}{2}|\dom(\eta)|\}\neq\emptyset.\]
Then $(z^*,\Delta^*)\in\K$ and $n(z^*,\Delta^*)\geq \lfloor\frac{1}{2}n(z,
\Delta)\rfloor$. 
\item Suppose that $(z_0,\Delta_0),\ldots,(z_n,\Delta_n)\in\K$ are such that
the sets $z_k$ (for $k\leq n$) are pairwise disjoint. Let $z=\bigcup\limits_{
k\leq n}z_k$, $\Delta=\bigcup\limits_{k\leq n}\Delta_k$. Then $(z,\Delta)\in
\K$ and $n(z,\Delta)\geq\min\{n(z_k,\Delta_k):k\leq n\}$.
\item Assume $(z,\Delta_0),(z,\Delta_1)\in\K$, $\Delta=\Delta_0\cup\Delta_1$. 
Then $(z,\Delta)\in\K$ and $n(z,\Delta)\geq\min\{\lfloor\frac{1}{2}n(z,
\Delta_0)\rfloor,\lfloor\frac{1}{2}n(z,\Delta_1)\rfloor\}$.
\end{enumerate}
\end{claim}

\noindent{\em Proof of the claim:}\quad 1)\ \ \ It should be clear that
$(z^*,\Delta^*)\in\K$. Suppose that $\Delta'\subseteq\Delta*$. For each $\nu
\in\Delta'$ fix $\eta_\nu\in\Delta$ such that
\[\nu=\eta \rest z^*\quad\mbox{ and }\quad|\dom(\eta_\nu)\cap z^*|\geq
\frac{1}{2}|\dom(\eta_\nu)|.\]
Loot at $\Delta^+\stackrel {\rm  def}{=}\{\eta_\nu:\nu\in\Delta'\}$. By the
definition of $n(z,\Delta)$ we find $\Delta^{++}\subseteq\Delta^+$ such that
elements of $\Delta^{++}$ have pairwise disjoint domains and
\[|\bigcup_{\eta_\nu\in\Delta^{++}}\dom(\eta_\nu)|\geq (z,\Delta)\cdot
|\Delta^+|.\]
So now look at $\Delta''\stackrel {\rm def}{=}\{\nu\in\Delta':\eta_\nu\in
\Delta^{++}\}$. Clearly, elements of $\Delta''$ have pairwise disjoint
domains and 
\[|\bigcup_{\nu\in\Delta''}\dom(\nu)|\geq\frac{1}{2}|\bigcup_{\eta_\nu\in
\Delta^{++}}\dom(\eta_\nu)|\geq\frac{n(z,\Delta)}{2}\cdot |\Delta^+|\geq
\lfloor\frac{1}{2}n(z,\Delta)\rfloor\cdot |\Delta^+|.\]

\medskip
\noindent 2)\ \ \ Suppose that $\Delta'\subseteq\Delta$. For $k\leq n$ let
$\Delta'_k=\Delta'\cap\Delta_k$ and choose $\Delta''_k\subseteq\Delta'_k$
such that the sets $\dom(\eta)$ (for $\eta\in\Delta''_k$) are pairwise
disjoint and 
\[|\bigcup_{\eta\in\Delta''_k}\dom(\eta)|\geq n(z_k,\Delta_k)\cdot
|\Delta'_k|.\]
Let $\Delta''=\bigcup\limits_{k\leq n}\Delta''_k$. Clearly the elements of
$\Delta''$ have pairwise disjoint domains. Moreover,
\[\begin{array}{l}
|\bigcup\limits_{\eta\in\Delta''}\dom(\eta)|=\sum\limits_{k\leq n}
\sum\limits_{\eta\in\Delta''_k}|\dom(\eta)|\geq \sum\limits_{k\leq n}n(z_k,
\Delta_k)\cdot |\Delta'_k|\\
\geq |\Delta'|\cdot \min\{n(z_k,\Delta_k): k\leq m\}.
\end{array}\]

\medskip
\noindent 3)\ \ \ Let $\Delta'\subseteq\Delta$. Let $\ell<2$ be such that
$|\Delta'\cap\Delta|\geq\frac{1}{2}|\Delta'|$. Now we may choose $\Delta''
\subseteq\Delta'\cap\Delta_\ell$ such that all members of $\Delta''$ have
pairwise disjoint domains and
\[|\bigcup_{\eta\in\Delta''}\dom(\eta)|\geq n(z,\Delta_\ell)\cdot
|\Delta'\cap\Delta_\ell|\geq \lfloor\frac{1}{2}n(z,\Delta_\ell)\rfloor\cdot
|\Delta'|.\]

\medskip
Now we are ready to define the triple $\trojka$. A semi--creature
$t\in\SCR[\bH]$ is taken to $K_{\ref{main}}$ if
\begin{itemize}
\item $\dis[t]=(z_t,\Delta_t)$, where either $\emptyset\neq z_t\in\fsuo$,
$\Delta_t=\emptyset$ or $(z_t,\Delta_t)\in\K$ and $n(z_t,\Delta_t)\geq 1$,
\item $\dom[t]=z_t$,
\item if $(z_t,\Delta_t)\in\K$ then 
\[\nor[t]=\log_8(n(z_t,\Delta_t))\quad\mbox{ and }\quad\sval[t]=v(z_t,
\Delta_t),\]
\item if $\Delta_t=\emptyset$ then $\nor[t]=\infty$ and $\sval[t]={\cal
X}^{\textstyle z_t}$.
\end{itemize}
For semi--creatures $t_0,\ldots,t_n\in K_{\ref{main}}$ with disjoint domains,
$\Sigma_{\ref{main}}(t_0,\ldots,t_n)$ consists of all $t\in K_{\ref{main}}$
such that
\[z_t=\bigcup_{\ell\leq n}z_{t_\ell}\quad\mbox{ and }\quad\Delta_t\supseteq
\bigcup_{\ell\leq n}\Delta_{t_\ell}\]
(so in particular $\bdom[t]=\bigcup\limits_{\ell\leq n}\bdom[t_\ell]$ and
$\sval[t]\subseteq\{v\in {\cal X}^{\textstyle \bdom[t]}: (\forall \ell\leq
n)(v \rest
\bdom[t_\ell]\in\sval[t_\ell])\}$). It should be clear that
$\Sigma_{\ref{main}}$ is a semi--composition operation on $K_{\ref{main}}$
(i.e.~the demands (a)--(c) of \ref{semi}(3) are satisfied). Next, for a
semi--creature $t\in K_{\ref{main}}$ we define $\Sigma^\bot_{\ref{main}}(t)$
as follows. It consists of all sets $\{s_0,\ldots,s_n\}\subseteq
K_{\ref{main}}$ such that
\[z_t=\bigcup_{\ell\leq n}z_{s_\ell}\quad\mbox{ and }\quad(\forall\eta\in
\Delta_t)(\exists\ell\leq n)(\eta \rest 
z_{s_\ell}\in\Delta_{s_\ell})\]
(note that this implies that at least one $\Delta_{s_\ell}$ is non-empty,
provided $\Delta_t\neq\emptyset$). Again, $\Sigma^\bot_{\ref{main}}$ is a
semi--decomposition operation on $K_{\ref{main}}$ (i.e.~clauses
(a)$^\bot$--(c)$^\bot$ of \ref{semi}(4) hold). Consequently, $\trojka$ is a
semi--creating triple for $\bH$, and plainly it is directly $+$--invariant and
directly permutation--invariant. It follows from \ref{cl2} that the forcing
notion $\forsing$ is not trivial. To show that the triple $\trojka$ has the
cutting property suppose that $t\in K_{\ref{main}}$, $\nor[t]>1$ and
$\emptyset\neq z
\varsubsetneq
\bdom[t]$. Let 
\[\begin{array}{l}
\Delta_0=\{\eta \rest 
z: \eta\in\Delta_t\ \&\ |\dom(\eta)\cap z|\geq\frac{1}{2}
|\dom(\eta)|\},\\
\Delta_1=\{\eta \rest 
(\bdom[t]\setminus z):\eta\in\Delta_t\ \&\ |\dom(\eta)
\setminus z|\geq\frac{1}{2}|\dom(\eta)|\}
\end{array}\]
and let $s_0,s_1\in K_{\ref{main}}$ be such that $\dis[s_0]=(z,\Delta_0)$,
$\dis[s_1]=(\bdom[t]\setminus z,\Delta_1)$. Easily $\{s_0,s_1\}\in
\Sigma^\bot_{\ref{main}}$ is as required in \ref{cccdef}(3) (remember
\ref{cl3}(1)). Similarly, it follows from \ref{cl3}(3) that $\trojka$ is
linked, and we immediately conclude from \ref{cl3}(2) that it is semi--gluing.


\medskip
\begin{conclusion} 
\label{conmain}
Let $({\cal X},+)$ be a finite Abelian group. Let $\trojka$ be the
semi--creating pair constructed in \ref{main} for $\bH(i)={\cal X}$ ($i\in 
\omega$). Then $\idmain$ is:
\begin{quotation}
\noindent a Borel ccc index--invariant translation--invariant
$\sigma$--ideal of subsets of ${\cal X}^{\textstyle\omega}$ which is
neither the meager ideal nor the null ideal nor their intersection. 
\end{quotation}
\end{conclusion}

\medskip
\begin{pf} 
\quad By \ref{invthm} and \ref{idccc} + \ref{cccthm} we know that $\idmain$
is a Borel ccc permutation--invariant translation--invariant $\sigma$--ideal
of subsets of ${\cal X}^{\textstyle\omega}$. Still we have to show the
following claim.

\begin{claim}
\label{cl5}
$\idmain$ is index--invariant. 
\end{claim}

\noindent{\em Proof of the claim:}\quad Let $\pi:\omega\stackrel{1-1}{
\longrightarrow}\omega$ be an embedding. Let
\[\begin{array}{r}
\q_\pi\stackrel{\rm def}{=}\{p\in\forsing: (\forall i\in\omega)(\bdom[t^p_i]
\cap\rng(\pi)=\emptyset\ \ \\
\mbox{or }\ \bdom[t^p_i]\subseteq\rng(\pi))\}.
  \end{array}\]
Like in \ref{cl1} one shows that $\q_\pi$ is a dense subset of $\forsing$
(remember that if $t_0,\ldots,t_n\in K_{\ref{main}}$ have disjoint domains and
norms $\infty$ then there is $t\in\Sigma_{\ref{main}}(t_0,\ldots,t_n)$ with
$\nor[t]=\log_8(n+1)$). Define a mapping $f_\pi:\q_\pi\longrightarrow\forsing$
by 
\[f_\pi(p)=(w^p\comp\pi,\langle\pi^{-1}(t^p_i):i\in\omega,\ \bdom(t^p_i)
\subseteq\rng(\pi)\rangle)\quad\mbox{ for } p\in\q_\pi.\]
(It should be clear that $f_\pi(p)$ is a condition in $\forsing$). We claim
that $f_\pi$ is a projection, and for this we have to show that
\begin{enumerate}
\item[$(\alpha)$] if $p,q\in\q_\pi$, $p\leq q$ then $f_\pi(p)\leq f_\pi(q)$,
and 
\item[$(\beta)$] if $p\in\q_\pi$, $r\in\forsing$ are such that $f_\pi(p)\leq
r$ then there is $q\in\q_\pi$ such that $p\leq q$ and $r\leq f_\pi(q)$.
\end{enumerate}
So suppose that $p,q\in\q_\pi$, $p\leq q$. Thus the condition $q$ can be
constructed from $p$ by repeating finitely many times the operations described
in \ref{forcing}: deciding the value, applying $\Sigma_{\ref{main}}$ and
applying $\Sigma^\bot_{\ref{main}}$. Now we would like to apply the same
procedures but restricting all creatures involved to $\rng(\pi)$ (remember
$p,q\in\q_\pi$). However it requires some extra care: in the original
procedure we may build at some moment a semi--creature $t$ such that
$\bdom[t]\cap\rng(\pi)\neq\emptyset\neq\bdom[t]\setminus\rng(\pi)$. What
should be the result $s$ of restricting $t$ to $\rng(\pi)$? We should take
$z_s=z_t\cap\rng(\pi)$, but we cannot take just $\Delta_s=\{\eta \rest
z_s:\eta\in\Delta_t\ \&\ \dom(\eta)\cap z_s\neq\emptyset\}$, as later (in some
application of $\Sigma^\bot_{\ref{main}}$) some functions $\eta\in\Delta_t$
may be restricted to sets disjoint from $\rng(\pi)$ (look at the definition of
$\Sigma^\bot_{\ref{main}}$). But we know where we are going to finish: the
active elements are $\nu\in\Delta_{t^q_i}$ for those $i$ that $\bdom[t^q_i]
\subseteq\rng(\pi)$. So, whenever in our procedure of building $q$ from $p$ we
create a semi--creature $t\in K_{\ref{main}}$ such that $\bdom[t]\cap\rng(\pi)
\neq\emptyset$, we replace it by $s$ such that $z_s=z_t\cap\rng(\pi)$ and
$\Delta_s\subseteq\{\eta\rest z_s:\eta\in\Delta_t\ \&\ (\exists\nu\in
\bigcup\limits_{i\in\omega} \Delta_{t^q_i})(\nu\subseteq\eta\rest z_s)\}$
is defined as follows. Let $\eta\in\Delta_t$. We ask if there is a
sequence of $\langle\eta_i,t_i:i\leq k\rangle$ such that $\eta_0=\eta$,
$\eta_i\in\Delta_{t_i}$, $\eta_i\supseteq\eta_{i+1}$, $\dom(\eta_k)\subseteq 
\rng(\pi)$, $t_0=t$, $t_k=t^q_j$ (for some $j$), $t_0,\ldots, t_k\in
K_{\ref{main}}$ appear at the successive levels of the construction (of
$q$ from $p$) after the one we consider. If the answer is positive we take
$\eta\rest z_s$ to $\Delta_s$, otherwise not. 

We ignore any $t$ with $\dom[t]\cap\rng(\pi)=\emptyset$ and whenever we apply
deciding the values we replace the respective sequence $w$ (extending $w^p$)
by $w \rest
\rng(\pi)$. This procedure results in transforming the sequence
\[(w^p \rest
\rng(\pi),t^p_i:\bdom[t^p_i]\subseteq\rng(\pi))\]
by successive applications of legal operations (and getting results with
domains included in $\rng(\pi)$). The final sequence is
$(w^q,t^q_i:\bdom[t^q_i]\subseteq\rng(\pi))$. Now we use $\pi^{-1}$ to
carry this procedure to $\forsing$ and we get a witness for $f_\pi (p)\leq
f_\pi(q)$.

To show clause $(\beta)$ suppose that $p\in\q_\pi$, $r\in\forsing$ are such
that $f_\pi(p)\leq r$. Let $w^q=(w^r\comp\pi^{-1})\cup (w^p\setminus\rng(
\pi))$, and $t^q_i=t^p_i$ whenever $\bdom[t^p_i]\cap\rng(\pi)=\emptyset$, and 
\[\{t^q_i:i\in\omega,\ \bdom[t^p_i]\subseteq\rng(\pi)\}=\{\pi(t^r_j): j\in
\omega\}.\]
It should be clear that $q=(w^q,t^q_0,t^q_1,\ldots)\in\q_\pi$ is stronger than
$p$ and $f_\pi(q)=r$.

Thus the mapping $f_\pi$ induces a complete embedding $g^*_\pi$ of the
complete Boolean algebra ${\rm BA}(\forsing)$ determined by $\forsing$ into
${\rm BA}(\q_\pi)$ (the last is, of course, isomorphic to ${\rm BA}(\forsing
)$). Assume that $A\in\idmain$ is a Borel set. We want to show that $\pi_*(A)
\in\idmain$. So suppose otherwise. Then we find a condition $p\in\q_\pi$
such that $p \forces\dot{W}\in\pi_*(A)$ or, in other words, $p \forces
\dot{W}\comp\pi\in A$. But $\dot{W}\comp\pi$ is the image of $\dot{W}$
under the embedding $g^*_\pi$, so $f_\pi(p)\forces\dot{W}\in A$, a
contradiction. Thus the claim is proved.

To finish the proof of the conclusion we have to check that the ideal
$\idmain$ is really new. For this it is enough to show the following.

\begin{claim}
\label{cl6}
\ 

\begin{enumerate}
\item For each $p\in\forsing$, $\POS(p)$ is a nowhere dense set.
\item For every $p\in\forsing$ there is $q\geq p$ such that $\POS(q)$ is a
null subset of ${\cal X}^\omega$.
\end{enumerate}
\end{claim}

\noindent{\em Proof of the claim:}\quad 1)\ \ \ Should be clear.

\medskip
\noindent 2)\ \ \ As $\trojka$ is semi--gluing we may assume that 
\[(\forall i\in\omega)(|\bdom[t^p_i]|>(3+i)\cdot |{\cal X}|^{2\cdot (3+i)}\
\mbox{ and }\ \nor[t^p_i]>5).\]
For each $i\in\omega$ choose a family $\Delta_i$ of partial functions with
disjoint domains contained in $\bdom[t^p_i]$ such that
\[|\Delta_i|=|{\cal X}|^{2\cdot (3+i)}\quad\mbox{ and }\quad\eta\in\Delta_i\
\Rightarrow |\eta|=3+i.\]
Let $s_i\in K_{\ref{main}}$ be such that
$\dis[s_i]=(z_{t^p_i},\Delta_i)$. Clearly $\nor[s_i]=\log_8(3+i)$, so
$q_0=(w^p,s_0,s_1,s_2,\ldots)\in\forsing$. Note that 
\[\Leb(\{v\in {\cal X}^\omega: (\forall\nu\in\Delta_i)(\nu\not\subseteq v)\})=
(1-\frac{1}{|{\cal X}|^{3+i}})^{|{\cal X}|^{2\cdot (3+i)}}\leq e^{-|{\cal
X}|^{3+i}}\]
(where $\Leb$ stands for the (product) Lebesgue measure on ${\cal X}^{
\textstyle \omega}$). Hence $\Leb(\POS(q_0))=0$. Since $\trojka$ is linked
(and $\bdom[s_i]=\bdom[t^p_i]$) we find a condition $q$ stronger than both $p$
and $q_0$. As $\POS(q)\subseteq\POS(q_0)$, we are done. 
\end{pf}

\bigskip
\stepcounter{section}
\subsection*{\quad 4. Making the example (perhaps) nicer}
The ideal $\idmain$ constructed in the previous section may be considered 
as a not so nice solution to \ref{pwa}. First note that there is no explicit
relation between the ideal and the sets $\POS(p)$. One would like to have here
that the family of all sets $\POS(p)$ is a basis of some ccc topology on
${\cal X}^{\textstyle \omega}$ and the ideal is the ideal of meager (with
respect to this topology) subsets of ${\cal X}^{\textstyle\omega}$, or at
least be able to apply Category Base approach of Morgan (see \cite{Mo77},
\cite{Mo90}). But it is not clear if we may get this in our example. Moreover,
the complexity of the incompatibility relation in $\forsing$ seems to be above
$\Sigma^1_1$, so the forcing notion is not really nice. However, we may modify
our example a little bit to get more nice properties of the forcing
notion. But the price to pay is that we go slightly out of the schema of {\em
norms on possibilities}.  

\medskip
\begin{definition}
\label{newexam}
\ 

\begin{enumerate}
\item Let $\N$ consists of all sequences $\bar{n}=\langle n^0_m,n^1_m:m<
\omega\rangle$ such that
\begin{itemize}
\item $4<n^0_m\leq n^1_m<n^0_{m+1}<\omega$  for each $m\in\omega$, and
\item $2^{2(m^*+2)}\cdot\sum\limits_{m<m^*}n^0_m\cdot n^1_m<n^0_{m^*}$ for
each $m^*\in\omega$, and 
\item $\lim\limits_{m\to\infty}(n^1_m)^{\frac{1}{2\cdot n^0_m}} =\infty$. 
\end{itemize}
\item Let $\bar{n}=\langle n^0_m,n^1_m:m<\omega\rangle\in\N$. We define a
forcing notion $\qnew$ (where $\bH:\omega\longrightarrow{\cal H}(\aleph_1)$)
as follows.

\noindent{\bf Conditions} are sequences $p=(w,\sigma_0,\sigma_1,\sigma_2,
\ldots)$ such that
\begin{enumerate}
\item[(a)] $w,\sigma_j$ (for $j\in\omega$) are finite functions with pairwise
disjoint domains, $w\in\prod\limits_{i\in\dom(w)}\bH(i)$, $\sigma_j\in
\prod\limits_{i\in\dom(\sigma_j)}\bH(i)$, 
\item[(b)] for some $m^*=m^*(p)<\omega$ there is a partition $\langle V^p_m:
m^*\leq m<\omega\rangle$ of $\omega$ such that for each $m\geq m^*$
\[|V^p_m|\leq n^1_m\cdot 2^{m^*}\quad\mbox{ and }\quad (\forall j\in V^p_m)(
|\dom(\sigma_j)|\geq\frac{n^0_m}{2^{m^*}}).\]
\end{enumerate}
For a condition $p=(w,\sigma_0,\sigma_1,\sigma_2,\ldots)$ we let
\[\POS(p)=\{\eta\in\prod_{i\in\omega}\bH(i): w\subseteq\eta\ \&\ (\forall
j\in\omega)(\sigma_j\not\subseteq\eta)\}.\]

\noindent{\bf The order} is given by:\qquad $p\leq q$ \quad if and only
if\quad $\POS(q)\subseteq\POS(p)$.
\item We will keep the convention that a condition $p\in\qnew$ is
$(w^p,\sigma ^p_0, \linebreak \sigma ^p_1,\ldots )$. 
\end{enumerate}
\end{definition}

\medskip
\begin{proposition}
\label{trivial}
Let $\bar{n}\in\N$.
\begin{enumerate}
\item $\qnew$ is a forcing notion.
\item For $p,q\in\qnew$, $p\leq q$\quad if and only if\quad $w^p\subseteq w^q$
and for each $i\in\omega$
\[(\exists j\in\omega)(\sigma^q_j\subseteq\sigma^p_i) 
\mbox{ or } 
(\exists m\in\dom(w^q)\cap\dom(\sigma^p_i))(w^q(m)\neq\sigma^p_i(m)).\]
\end{enumerate}
\end{proposition}

\begin{pf}
\quad  2)\ \ \ Assume that $\POS(q)\subseteq\POS(p)$.\\
First note that necessarily $w^q\cup w^p$ is a function. Suppose that
$k\in\dom(w^p)\setminus\dom(w^q)$. If $k\notin\bigcup\limits_{i\in\omega}
\dom(\sigma^q_i)$ then we may easily construct a function in $\POS(q)\setminus
\POS(p)$. So for some $i$ we have $k\in\dom(\sigma^q_i)$. Take $\ell\in\dom(
\sigma^q_i)\setminus\{k\}$ and build a function $\eta\in\POS(q)$ such that
$\eta(\ell)\neq\sigma^q_i(\ell)$, $\eta(k)\neq w^q(k)$. Then $\eta\notin
\POS(p)$, a contradiction showing that $w^p\subseteq w^q$.\\
Suppose now that $i\in\omega$ is such that for no $j\in\omega$, $\sigma^q_j
\subseteq\sigma^p_i$. Then for each $j\in\omega$ such that $\dom(\sigma^q_j)
\cap\dom(\sigma^p_i)\neq\emptyset$, either there is $k\in\dom(\sigma^q_j)
\cap\dom(\sigma^p_i)$ such that $\sigma^q_j(k)\neq\sigma^p_i(k)$ or
$\dom(\sigma^q_j)\setminus\dom(\sigma^p_i)\neq\emptyset$. Thus we may build a
partial function $\eta\in\prod\limits_{i\in\dom(\eta)}\bH(i)$ with
$\dom(\eta)=\bigcup\limits_{j\in\omega}\dom(\sigma^q_j)$ and such that
\[\eta \rest
\dom(\sigma^p_i)\subseteq\sigma^p_i\quad\mbox{ and }\quad(\forall
j\in\omega)(\sigma^q_j\not\subseteq\eta).\]
As any member of $\prod\limits_{i\in\omega}\bH(i)$ extending $\eta\cup w^q$ is
in $\POS(q)$ (and thus in $\POS(p)$) we conclude that the functions $w^q$ and
$\sigma^p_i$ are incompatible, i.e.
\[(\exists m\in\dom(w^q)\cap\dom(\sigma^p_i))(w^q(m)\neq\sigma^p_i(m)).\]
The converse implication is even easier. 
\end{pf}

\medskip
Let us recall that a partial order $(\q,\leq_\q)$ is Souslin (Borel,
respectively) if $\q$, $\leq_\q$ and the incompatibility relation $\bot_{\q}$
are $\Sigma^1_1$ (Borel, respectively) subsets of $\r$ and $\r\times\r$. On
Souslin forcing notions and their applications see Judah Shelah
\cite{JdSh:292}, Goldstern Judah \cite{GoJu} and Judah Ros{\l}anowski Shelah
\cite{JRSh:373} (the results of these three and many other papers on the topic
are presented in Bartoszy\'nski Judah \cite{BaJu95} too). A new systematic
treatment of definable forcing notions is presented in a forthcoming paper
\cite{Sh:630} (several results there are applicable to forcing notions
defined in this paper).

\medskip
\begin{theorem}
\label{niceforcing}
Suppose that $\bar{n}\in\N$ and $\bH:\omega\longrightarrow{\cal H}(\aleph_1)$.
\begin{enumerate}
\item The forcing notion $\qnew$ is $\sigma$-$*$--linked.
\item Conditions $p_0,p_1\in\qnew$ are compatible (in $\qnew$)\quad if and
only if\\
$\POS(p_0)\cap\POS(p_1)\neq\emptyset$.
\item The forcing notion $\qnew$ is Souslin ccc.
\end{enumerate}
\end{theorem}

\medskip
\begin{pf}
\quad 1)\ \ \ For a condition $p\in\qnew$ let $m^*(p)$ and $\langle V^p_m:
m^*(p)\leq m<\omega\rangle$ be given by \ref{newexam}(2b).\\
Fix $n\in\omega$.\\
For $m^*\in\omega$, $w\in\prod\limits_{i\in\dom(w)}\bH(i)$, $\dom(w)\in\fsuo$
and sequences
\[\bar{V}=\langle V_m: m^*\leq m<m^*+n+2\rangle\subseteq\fsuo\quad\mbox{ and
}\quad \bar{\sigma}=\langle\sigma_j: j\in\bigcup\limits_{m=m^*}^{m^*+n+1} V_m
\rangle\]
we let 
\[\begin{array}{ll}
A^{\bar{V},\bar{\sigma}}_{m^*,w}=\{p\in\qnew: &m^*(p)=m^*\ \&\ (\forall m\in
[m^*,m^*{+}n{+}2))(V^p_m=V_m)\ \&\\
\ &(\forall j\in\bigcup\limits_{m=m^*}^{m^*+n+1}V_m)(\sigma^p_j=\sigma_j)\ \&\
w^p=w\}.
  \end{array}\]
We want to show that any $n+1$ members of $A^{\bar{V},\bar{\sigma}}_{m^*,w}$
have a common upper  bound in $\qnew$. For this we will need the following
two technical observations.

\begin{claim}
\label{cl7}
Suppose that $p_0,\ldots,p_n\in\qnew$, $\max\{m^*(p_0),\ldots,m^*(p_n)\}\leq
m^*$ and $U\subseteq\bigcup\{V^{p_\ell}_m\times\{\ell\}:\ell\leq n,\ m^*+n+2
\leq m\}$ is finite. Then there are pairwise disjoint sets $u_{j,\ell}$ for
$(j,\ell)\in U$ such that  
\[u_{j,\ell}\subseteq\dom(\sigma^{p_\ell}_j)\ \mbox{ and }\ |u_{j,\ell}|\geq
\frac{n^0_m}{2^{m^*+n+1}},\]
where $m$ is such that $j\in V^{p_\ell}_m$.
\end{claim}

\noindent{\em Proof of the claim:}\quad Let $y(j,\ell)=\lfloor\frac{n^0_m}{2^{
m^*+n+1}}\rfloor+1$ (for $(j,\ell)\in U$ and $m$ such that $j\in
V^{p_\ell}_m$). For $(j,\ell)\in U$ and $y<y(j,\ell)$ let $a^y_{j,\ell}=\dom(
\sigma^{p_\ell}_j)$. We want to apply Hall's theorem to the system $\langle
a^y_{j,\ell}: (j,\ell)\in U\ \&\ y<y(j,\ell)\rangle$. So suppose that
$\A\subseteq\{(j,\ell,y):(j,\ell)\in U\ \&\ y<y(j,\ell)\}$. For some $\ell^*
\leq n$ we have $|\{(j,\ell,y)\in\A:\ell=\ell^*\}|\geq\frac{|\A|}{n+1}$. Now,
remembering that
\[|a^y_{j,\ell}|\geq\frac{n^0_m}{2^{m^*(p_\ell)}}\geq \frac{n^0_m}{2^{m^*}}>
y(j,\ell)\cdot 2^n,\]
we easily conclude that
\[|\bigcup\{a^y_{j,\ell^*}: (j,\ell^*,y)\in\A\}|\geq|\{(j,\ell,y)\in\A:
\ell=\ell^*\}|\cdot 2^n\geq\frac{|\A|}{n+1}\cdot 2^n\geq |\A|.\]
Consequently we may apply the Hall theorem (see \cite{Ha35}) and choose a
system $\langle k^y_{j,\ell}:(j,\ell)\in U\ \&\ y<y(j,\ell)\rangle$ of
distinct representatives for $\langle a^y_{j,\ell}:(j,\ell)\in U,\ y<y(j,\ell)
\rangle$ (so $k^y_{j,\ell}\in\dom(\sigma^{p_\ell}_j)$). Now let $u_{j,\ell}=
\{k^y_{j,\ell}: y<y(j,\ell)\}$ (for $(j,\ell)\in U$). It should be clear that
these sets are as required.

\begin{claim}
\label{cl8}
Suppose $p_0,\ldots,p_n\in\qnew$, $\max\{m^*(p_0),\ldots,m^*(p_n)\}\leq m^*$. 
There is a sequence $\langle u_{j,\ell}: \ell\leq n,\ j\in\bigcup\{ 
V^{p_\ell}_m:m\geq m^*+n+2\}\rangle$ of pairwise disjoint sets such that  
\[u_{j,\ell}\subseteq\dom(\sigma^{p_\ell}_j)\ \mbox{ and }\ |u_{j,\ell}|\geq
\frac{n^0_m}{2^{m^*+n+1}},\]
where $m$ is such that $j\in V^{p_\ell}_m$.
\end{claim}

\noindent{\em Proof of the claim:}\quad By \ref{cl7} we know that for each
finite $U\subseteq\bigcup\{V^{p_\ell}_m\times\{\ell\}:\ell\leq n,\ m^*+n+2
\leq m\}$ we can find a sequence $\langle u_{j,\ell}: (j,\ell)\in U\rangle$
with the respective properties. Of course, for each $U^*\subseteq U$ the
restricted sequence $\langle u_{j,\ell}: (j,\ell)\in U^*\rangle$ will have
those properties too. Moreover, for each finite $U$ the number of all possible
sequences is finite. Consequently we may use K\"onig lemma and conclude that
there exists $\langle u_{j,\ell}: \ell\leq n,\ j\in\bigcup\{V^{p_\ell}_m:m\geq
m^*+n+2\}\rangle$ as desired.

\begin{claim}
\label{cl9}
Suppose $p_0,\ldots,p_n\in A^{\bar{V},\bar{\sigma}}_{m^*,w}$. Then the
conditions $p_0,\ldots,p_n$ have a common upper bound in $\qnew$.
\end{claim}

\noindent{\em Proof of the claim:}\quad Using \ref{cl8} choose a sequence
\[\langle u_{j,\ell}: \ell\leq n\ \&\ j\in\bigcup\{V^{p_\ell}_m:m\geq m^*+n+2
\}\rangle\]
of pairwise disjoint sets such that $u_{j,\ell}\subseteq\dom(\sigma^{
p_\ell}_j)$ and $|u_{j,\ell}|\geq\frac{n^0_m}{2^{m^*+n+1}}$, where $m$ is such
that $j\in V^{p_\ell}_m$. Let $\sigma^q_{j,\ell}=\sigma^{p_\ell}_j \rest
u_{j,\ell}$ and let $w^q$ be a finite function such that $w^q\in\prod
\limits_{i\in\dom(w^q)}\bH(i)$,
\[\dom(w^q)=\dom(w)\cup\bigcup\{\dom(\sigma_j): j\in\bigcup_{m=m^*}^{m^*+n+1}
V_m\}\]
and $w\subseteq w^q$ and if $k\in\dom(\sigma_j)$, $j\in V_m$, $m^*\leq m<m^*+n
+2$ then $\sigma_j(k)\neq w^q(k)$. Look at the sequence 
\[q=(w^q,\sigma^q_{j,\ell}:\ell\leq n,\ j\in\bigcup\{V^{p_\ell}_m: m^*+n+2\leq
m\}).\]
It is a condition in $\qnew$ with $m^*(q)=m^*+n+2$ and $V^q_m=\{(j,\ell):\ell
\leq n\ \&\ j\in V^{p_\ell}_m\}$ witnessing the clause \ref{newexam}(2b). (To
be strict, one should re-enumerate all $\sigma^q_{j,\ell}$ to have single
indexes, but that is not a problem.) Immediately by the choice of the
$\sigma^q_{j,\ell}$'s and $w^q$ one concludes
\[\POS(q)\subseteq\POS(p_0)\cap\ldots\cap\POS(p_n),\]
finishing the proof of the claim.

\medskip
Since there are countably many possibilities for indexes
$(\bar{V},\bar{\sigma},m^*,w)$ (in $A^{\bar{V},\bar{\sigma}}_{m^*,w}$), the
first part of the theorem follows from \ref{cl9}.

\medskip
\noindent 2)\ \ \ Suppose $\eta\in\POS(p_0)\cap\POS(p_1)$. Let $m^*\geq
m^*(p_0)+m^*(p_1)+3$ be such that $2^{m^*+1}>|\dom(w^{p_0})\cup\dom(w^{p_1}
)|$. For $\ell<2$, $m\in [m^*(p_\ell),m^*)$ and $j\in V^{p_\ell}_m$ choose an
integer $k^\ell_j\in\dom(\sigma^{p_\ell}_j)$ such that $\eta(k^\ell_j)\neq
\sigma^{p_\ell}_j(k^\ell_j)$. Let 
\[a=\dom(w^{p_0})\cup\dom(w^{p_1})\cup \{k^\ell_j:\ell<2,\ j\in
\bigcup_{m=m^*(p_\ell)}^{m^*-1}V^{p_\ell}_m\}.\]
Plainly,
\[|a|\leq |\dom(w^{p_0})\cup\dom(w^{p_1})|+2^{m^*+1}\cdot\sum_{m<m^*}
n^1_m<\frac{n^0_{m^*}}{2^{m^*}}\]
(remember \ref{newexam}(1)). Now for $\ell<2$, $m\geq m^*$ and $j\in
V^{p_\ell}_m$ let $\sigma^*_{\ell,j}=\sigma^{p_\ell}_j \rest
(\dom(\sigma^{p_\ell}_j)\setminus a)$. Note that (for relevant $\ell,j,m$)
\[|\dom(\sigma^*_{\ell,j})|>\frac{n^0_m}{2^{m^*(p_\ell)}}-\frac{n^0_{m^*}}{
2^{m^*}}>4\frac{n^0_m}{2^{m^*}}.\]
Like in \ref{cl9}, use \ref{cl8} to get a sequence 
\[(\sigma^q_{j,\ell}:\ell<2,\ j\in\bigcup\{V^{p_\ell}_m: m^*\leq m\})\]
such that $\sigma^q_{j,\ell}\subseteq\sigma^*_{\ell,j}$, $|\dom(\sigma^q_{j,
\ell})|\geq\frac{n^0_m}{2^{m^*}}$. Next we let $w^q=\eta \rest 
a$ and as in
\ref{cl8} we see that
\[q=(w^q,\sigma^q_{j,\ell}:\ell<2,\ j\in\bigcup\{V^{p_\ell}_m: m^*\leq
m\})\in\qnew\]
is a condition stronger than both $p_0$ and $p_1$. (But note that we cannot
claim that $\eta\in\POS(q)$.).

\medskip
\noindent 3)\ \ \ Let ${\cal Y}$ be the space of all sequences $(w,\sigma_0,
\sigma_1,\ldots)$ of finite functions such that $w\in\prod\limits_{i\in
\dom(w)}\bH(i)$,
$\sigma_j\in\prod\limits_{i\in\dom(\sigma_j)}\bH(i)$. The space ${\cal Y}$,
equipped with the product topology of discrete spaces, is a Polish space. Now,
$\qnew$ is a $\Sigma^1_1$ subset of ${\cal Y}$, as to express that
$(w,\sigma_0,\sigma_1,\ldots)\in\qnew$ we need to say that ``{\em there is} a
partition $\langle V_m: m^*\leq m<\omega\rangle$ of $\omega$ as in
\ref{newexam}(2b)'' and the rest of the demands is Borel. The relation
$\leq_{\qnew}$ is clearly $\Sigma^1_1$ if one uses \ref{trivial}(2) to express
it. Finally, to deal with the incompatibility relation look at the proof of
clause 2) above. To say that conditions $p_0,p_1$ are compatible we have to
say that, for $m^*$ as there, we can find points $k^\ell_j\in\dom(\sigma^{
p_\ell})$ (for $\ell<2$, $j$ as there) and a function $w^q$ such that
$w^{p_0}\cup w^{p_1}\subseteq w^q$ and $w^q(k^\ell_j)\neq\sigma^{p_\ell}_j
(k^\ell_j)$. 
\end{pf}

\medskip
\begin{remark}
\label{remark}
Note that one can easily choose a dense suborder $\q^*\subseteq\qnew$ such
that $\q^*$ is Borel ccc: take to $\q^*$ those conditions $p$ for which
\linebreak
$|\dom(\sigma_j)|= \lfloor\frac{n^0_m}{2^{m^*(p)}}\rfloor+1$ for $j\in V^p_m$,
$m\geq m^*(p)$. 
\end{remark}
 
\medskip
Before we introduce an ideal related to $\qnew$ let us recall some basic
notions of the Category Base technique.

\medskip
\begin{definition}
[Morgan, see \cite{Mo77}, \cite{Mo90}]
\ 

\begin{enumerate}
\item A family $\C$ of subsets of a space $\cal X$ is called {\em a category
base on $\cal X$} if
\begin{enumerate}
\item[(a)] ${\cal X}=\bigcup\C$,
\item[(b)] for every $\D\subseteq\C$ consisting of disjoint sets and such that
$0<|\D|<|\C|$ and for any $A\in\C$:

if  $(\exists B\in\C)(B\subseteq A\cap\bigcup\C)$ then $(\exists D\in\D)
(\exists B\in\C)(B\subseteq A\cap D)$, and

if $\neg(\exists B\in\C)(B\subseteq A\cap\bigcup\C)$ then $(\exists B\in\C)
(B\subseteq A\setminus \bigcup\D)$.
\end{enumerate}
\item Let $\C$ be a category base on $\cal X$. A set $X\subseteq{\cal X}$ is
called 
\begin{itemize}
\item {\em $\C$--singular} if $(\forall A\in\C)(\exists B\in\C)(B\subseteq
A\setminus X)$,
\item {\em $\C$--meager} if $X$ can be covered by a countable union of
$\C$--singular sets.
\end{itemize}
We say that a set $X\subseteq{\cal X}$ has {\em the $\C$--Baire property} if
for every $A\in\C$ there is $B\in\C$, $B\subseteq A$ such that either $B\cap
X$ or $B\setminus X$ is $\C$--meager.
\item For a category base $\C$ on $\cal X$, the family of all $\C$--meager
sets will be denoted by $\M_\C$ and the family of subsets of $\cal X$ with the
$\C$--Baire property will be called $\B_\C$.
\end{enumerate}
\end{definition}

\medskip
\begin{remark}
A category base, $\C$--meager sets and sets with $\C$--Baire property
generalize the notions of a topology, meager sets and sets with the Baire
property (with respect to the topology). Several results true in the
topological case remain true for the category base approach. In particular,
$\B_\C$ is a $\sigma$-field of subsets of $\cal X$ closed under the Souslin
operation $\A$, and $\M_\C$ is a $\sigma$-ideal of subsets of $\cal X$. For a
systematic study of the category base method we refer the reader to Morgan
\cite{Mo90}. 
\end{remark}

\medskip
\begin{conclusion}
\label{another}
Let $\bar{n}\in\N$ and $\bH:\omega\longrightarrow{\cal H}(\aleph_1)$. 
\begin{enumerate}
\item The family $\C^{\bH}_{\bar{n}}\stackrel{\rm def}{=}\{\POS(p): p\in
\qnew\}$ is a category base on $\prod\limits_{i\in\omega}\bH(i)$ such that no
member of $\C^{\bH}_{\bar{n}}$ is $\C^{\bH}_{\bar{n}}$--meager. 
\item Each $\Sigma^1_1$--subset of $\prod\limits_{i\in\omega}\bH(i)$ has the
$\C^{\bH}_{\bar{n}}$--Baire property. For every set $A\subseteq\prod\limits_{
i\in\omega}\bH(i)$ with the $\C^{\bH}_{\bar{n}}$--Baire property there is a
$\Pi^0_3$--subset $B$ of $A$ such that $A\setminus B$ is
$\C^{\bH}_{\bar{n}}$--meager. 
\item The $\sigma$-ideal $\M_{\C^{\bH}_{\bar{n}}}$ of
$\C^{\bH}_{\bar{n}}$--meager subsets of $\prod\limits_{i\in\omega}\bH(i)$ is a
Borel ccc $\sigma$-ideal (in fact, any member of $\M_{\C^{\bH}_{\bar{n}}}$ can
be covered by a $\Sigma^0_3$--set from $\M_{\C^{\bH}_{\bar{n}}}$). 
\item Suppose that $\bH(i)={\cal X}$ (for $i\in\omega$) and $\cal X$ is a
finite group. Then $\M_{\C^{\bH}_{\bar{n}}}$ is a Borel ccc
translation--invariant index--invariant $\sigma$-ideal of subsets of ${\cal
X}^{\textstyle \omega}$, which is neither the ideal of meager sets, nor the
ideal of null sets nor their intersection. Moreover, the formula {\em ``a real
$r$ a is a Borel code for a subset of ${\cal X}^\omega$ and the set $\# r$
coded by $r$ is in $\M_{\C^{\bH}_{\bar{n}}}$''} is $\Sigma^1_2$.
\end{enumerate}
\end{conclusion}

\medskip
\begin{pf}
\quad 1)--3)\ \ \ Should be clear if you remember \ref{niceforcing}(2,3). 
Note that each $\POS(p)$ is closed.

\medskip
\noindent 4)\ \ \ Plainly, the ideal $\M_{\C^{\bH}_{\bar{n}}}$ is
translation--invariant. To estimate the complexity of the formula ``$A\in
\M_{\C^{\bH}_{\bar{n}}}$'' use remark \ref{remark}. The remaining assertions
follow from the following two observations.

\begin{claim}
If $A\subseteq{\cal X}^{\textstyle \omega}$ is $\C^\bH_{\bar{n}}$--singular
and $\pi:\omega\stackrel{1-1}{\longrightarrow}\omega$ is an embedding then
$\pi_*(A)$ is $\C^\bH_{\bar{n}}$--singular.
\end{claim}

\noindent{\em Proof of the claim:}\quad Let $p\in\qnew$. Passing to a stronger
condition if needed we may assume that
\[(\forall j\in\omega)(\dom(\sigma^p_j)\subseteq\rng(\pi)\ \mbox{ or }\
\dom(\sigma^p_j)\cap\rng(\pi)).\]
Let $q=(w^p\comp\pi,\sigma^p_j\comp\pi:\ j\in\omega,\dom(\sigma^p_j)\subseteq
\rng(\pi))$. Plainly $q\in\qnew$, so we find $r\in\qnew$ such that $\POS(r)
\subseteq\POS(q)\setminus A$. Let $m^*=m^*(p)+m^*(r)+1$. Choose $w^{p^*}$ such
that 
\begin{enumerate}
\item $\dom(w^{p^*})=\dom(w^p)\cup\pi[\dom(w^r)]\cup\bigcup\{\dom(\sigma^p_j):
\dom(\sigma^p_j)\cap\rng(\pi)=\emptyset\ \&\ j\in\bigcup\limits_{m=m^*(p)}^{
m^*-1} V^p_m\}\cup\bigcup\{\pi[\dom(\sigma^r_i)]: i\in\bigcup\limits_{m=m^*(r)
}^{m^*-1} V^r_m\}$,
\item $w^p\rest 
(\omega\setminus\rng(\pi))\cup w^r\comp\pi^{-1}\subseteq
w^{p^*}$, and
\item $\sigma^p_j\not\subseteq w^{p^*}$, $\sigma^r_i\comp\pi^{-1}\not\subseteq
w^{p^*}$ (for $j,i$ as in 1) above).
\end{enumerate}
Next choose $\sigma^{p^*}_j$ such that $\langle\sigma^{p^*}_j:j<\omega\rangle$
enumerates the set
\[\{\sigma^p_j:\rng(\sigma^p_j)\cap\rng(\pi)=\emptyset\ \&\ j\in\bigcup_{m\geq
m^*}V^p_m\}\cup\{\sigma^r_j\comp\pi^{-1}:j\in\bigcup_{m\geq m^*}V^r_m\}.\] 
Plainly, $p^*=(w^{p^*},\sigma^{p^*}_0,\sigma^{p^*}_1,\ldots)\in\qnew$ and
$\POS(p^*)\subseteq\POS(p)\setminus \pi_*(A)$.

\begin{claim}
For each $p\in\qnew$ there is $q\in\qnew$ such that $\POS(q)\subseteq\POS(p)$
and $\POS(q)$ is nowhere dense and null.
\end{claim}

\noindent{\em Proof of the claim:}\quad Take $m^*>m^*(p)+5$ such that 
\[(\forall k\geq m^*)(n^1_k>|{\cal X}|^{2\cdot n^0_k})\]
and choose a sequence $\langle\sigma_{k,j}:k\geq m^*,j<n^1_k\rangle$ of finite
functions with pairwise disjoint domains and such that $|\dom(\sigma_{k,j})|=
n^0_k$ and 
\[\dom(\sigma_{k,j})\cap(\dom(w^p)\cup\bigcup\{\dom(\sigma^p_i): i\in
\bigcup_{m=m^*(p)}^{m^*-1}V^p_m\})=\emptyset.\]
Let $q=(w^p,\sigma_{k,j}: k\geq m^*,j<n^1_k)$. Easily $q\in\qnew$ and the
conditions $p,q$ are compatible (compare the proof of
\ref{niceforcing}(1)). Note that for each $k\geq m^*$
\[\begin{array}{l}
\Leb(\{\eta\in {\cal X}^{\textstyle \omega}: (\forall j<n^1_k)(\sigma_{k,j}
\not\subseteq\eta)\})=(1-\frac{1}{|{\cal X}|^{n^0_k}})^{n^1_k}<\\
(1-\frac{1}{|{\cal X}|^{n^0_k}})^{|{\cal X}|^{2\cdot n^0_k}}\leq e^{-|{\cal
X}|^{n^0_k}}.
  \end{array}\] 
Hence $\POS(q)$ is a nowhere dense null set and now we easily finish. 
\end{pf}

\medskip
\begin{remark}
There are several possible variants of the forcing notions $\qnew$, each of
them doing the job. These forcing notions will be presented in a subsequent
paper \cite{Sh:F224}, where we will systematically study forcing properties
of ccc partial orders build by the method of norms on possibilities. In
particular, we will show there that we may get different forcing properties of
$\qnew$ (for different $\bar{n}$), thus showing that the corresponding ideals
are not isomorphic. Let us note that the forcing notions which appeared in
the present paper are not $\sigma$--centered and do add Cohen reals. The proof
of these facts and more general statements will be presented in \cite{Sh:F224}.
\end{remark}


\bigskip
\bigskip
\noindent
{\small \uppercase {
Andrzej Ros{\l}anowski\hfill Saharon Shelah\\
Institute of Mathematics \hfill
Institute of Mathematics\\
The Hebrew University \hfill  The Hebrew University\\
of Jerusalem \hfill
of Jerusalem \\
 91904 Jerusalem, Israel  \hfill  91904 Jerusalem, Israel\\
{\lowercase {and\hfill and}}\\ 
Mathematical Institute  \hfill
Department of Mathematics\\
of Wroclaw University\hfill
Rutgers University\\
Pl. Grunwaldzki 2/4 \hfill New Brunswick, NJ 08854\\
50-384 Wroc{\l}aw, Poland \hfill
USA\\
roslanow@@math.huji.ac.il \hfill
shelah@@math.huji.ac.il\\
}}
\end{document}